
\input amstex \documentstyle {amsppt}
\magnification=1200
\hsize=15truecm
\baselineskip=12pt
\hoffset=0truecm
\predefine\lll{\l}
\predefine\ssn{\ss}

\redefine\cA{{\Cal A}}
\redefine\cB{{\Cal B}}
\redefine\cC{{\Cal C}}

\redefine\cF{{\Cal F}}
\redefine\cV{{\Cal V}}

\redefine\Pz{{\Cal P}_2}
\redefine\ss{\text{{\bf s}}}

\redefine\CC{\text{{\bf C}}}

\redefine\HH{{\Cal H}}

\redefine\ff{\varphi}

\redefine\cA{{\Cal A}}
\redefine\cB{{\Cal B}}
\redefine\cC{{\Cal C}}
\redefine\cF{{\Cal F}}

\redefine\HH{{\Cal H}}

\redefine\lb{\lbrack}
\redefine\rb{\rbrack}
\redefine\l{$\lbrack$}
\redefine\r{$\rbrack$}

\redefine\cV{{\Cal V}}
\redefine\CC{\text{{\bf C}}}
\redefine\RR{{\text{{\bf R}}}}
\redefine\NN{\text{{\bf N}}}

\redefine\sqp{\sqcup \hskip -0.9em \sqcap \hskip -0.95em +}

\redefine\ct{{\text{{\bf t}}}}
\redefine\fen{f_1\otimes\dots\otimes f_n}
\redefine\gen{g_1\otimes\dots\otimes g_n}
\redefine\gem{g_1\otimes\dots\otimes g_m}
\redefine\fzn{f_2\otimes\dots\otimes f_n}
\redefine\Pzr{\Pz(1,\dots,2r)}
\redefine\Pza{\Pz^{NC}(1,\dots,2r)}
\redefine\Pzaa{\Pz^{NC}}
\redefine\Pzn{\Pz(1,\dots,n)}
\redefine\PzN{\Pz(\infty)}
\redefine\oo{\omega}
\redefine\rt{\rho_\ct}
\redefine\rN{\rho^N}
\redefine\rte{\rho_{\ct_1}}
\redefine\rtz{\rho_{\ct_2}}
\redefine\rtQh{\rho_\ct^Q}
\redefine\rqe{\rho_{q_1}}
\redefine\rqz{\rho_{q_2}}
\redefine\rqcc{\rho_q\lb c^\sharp\dots c^\sharp\rb}
\redefine\cAn{\cA^0}

\redefine\cCn{\cC^0}
\redefine\odo{\otimes\dots\otimes}
\redefine\koi{\langle\{\oo_i\}\rangle}
\redefine\kci{\langle\{c_i\}\rangle}
\redefine\kbi{\langle\{b_i\}\rangle}
\document
\heading
{\bf INTERPOLATIONS BETWEEN} \\ {\bf BOSONIC AND FERMIONIC RELATIONS} \\
{\bf GIVEN BY} \\ {\bf GENERALIZED BROWNIAN MOTIONS} \\
\quad\\
{\bf Marek Bo\D zejko and Roland Speicher}\\
\quad\\ \quad\\
Instytut Matematyczny\\
Uniwersytet Wroc\lll awski\\ Plac Grunwaldzki 2/4\\
50-384 Wroc\lll aw\\ Poland\\
bozejko \@ math.uni.wroc.pl \\
\quad\\ and\\ \quad\\
Institut f\"ur Angewandte Mathematik\\ Universit\"at Heidelberg\\
Im Neuenheimer Feld 294\\ D-69120 Heidelberg\\
Federal Republic of Germany\\
L95 \@ vm.urz.uni-heidelberg.de
\endheading
\pagebreak
\noindent  running head: Generalized Brownian Motions
\newline
\vskip2cm\noindent
Please send all correspondence to:\newline
Roland Speicher:
Institut f\"ur Angewandte Mathematik, Universit\"at Heidelberg,
Im Neuenheimer Feld 294,
D-69120 Heidelberg, Federal Republic of Germany,\newline
e-mail: L95 \@ vm.urz.uni-heidelberg.de
\newline\pagebreak
\heading
{\bf Abstract}
\endheading
We present an interpolation between the bosonic and fermionic
relations. This interpolation is given by an object which we call
\lq generalized Brownian motion' and which is characterized by a
generalization of the pairing rule for the calculation of the moments
of bosonic and fermionic fields. We develop some basic theory for such
generalized Brownian motions and consider more closely one example,
which turns out to be intimately connected with Voiculescu's concept
of \lq free product'.\newline
\pagebreak
\baselineskip=12pt
\heading
{\bf 1. Introduction}
\endheading
We shall present here an interpolation between the bosonic and fermionic
relations. Such interpolations have attracted some attention in connection
with quantum groups, comp., e.g., \l Gre,Fiv,LPo\r.\par
As in \l BSp1\r, where we considered another example of such an
interpolation, our work is motivated from a probabalistic point of view.
Thus we are led to objects which we call \lq generalized Brownian motions'.
We shall now give a short survey on our probabalistic motivation and develop
some general theory on these generalized Brownian motions, whereas, in
Sect. 3, we shall give the construction of our special interpolation.
\par
In non-commutative (quantum) probability theory we are in search of
non-commutative generalizations of the classical probabilistic notions and
concepts. In particular, we are interested in generalizations of
processes with independent and stationary increments (\lq white noises'
or \lq Brownian motions') and the corresponding stochastic
integration theories (Ito-formulas). There are some general theories on
these objects \l AFQ,KPr,K\"um1\r,
but we believe that at the moment we are in need
of some more concrete examples of such white noises in order to get a
feeling for the typical properties and difficulties arising in this field,
before we can  hope to develop the general theory any further.\par
In \l BSp1,2\r\ we started the construction of such special examples
of Brownian motions. Here, we shall present another example.
All our examples belong to a more general class of Brownian motions,
which arise on one hand via some general central limit theorem relying
on the notion of \lq generalized independence' of K\"ummerer \l K\"um2\r\
and which on the other hand are motivated by the fact, that the
moments of classical Brownian motion can be calculated by pair partitions.
We can characterize our class of Brownian motions by a generalization of
this pairing in the following formal way: Let $t\to \omega(t)$ be our
non-commutative Brownian motion and consider its \lq increments'
$\omega(f)=\int f(t)d\omega(t)$ for $f\in L^2(\RR)$. Then we want to
define a state $\rt$ on the $*$-algebra generated by all $\omega(f)$
(with $\omega(f)^*=\omega(\bar f)$)
by the following \lq pairing prescription':
 $$
\rt\lb \omega(f_1)\dots \omega(f_n)\rb =
\cases 0,&\text{if $n$ odd}\\
\sum_{\cV=\{V_1,\dots,V_r\}}
\rho\lb V_1\rb\dots \rho\lb V_r\rb \cdot \ct(\cV),&
\text{if $n=2r$,}\endcases$$
where the sum runs over all 2-partitions (pairings) of $\{1,\dots,2r\}$
into sets
$V_1,$ $\dots,$ $V_r$, each consisting of two elements, and where
$\rho\lb V\rb$ denotes $<\bar f_i,f_j>$ for
$V=(i,j)$, and $\ct$ is some function on the set of all 2-partitions.
The specification of $\ct$ thus determines the concrete structure of
our Brownian motion. The main problem in this context is to decide
whether the linear functional $\rt$ is positive, i.e. whether it is
indeed a state. If this is the case, then we shall call $\ct$
\lq positive definite'. Note that this is a quite indirect definition,
but till now we have not been able to connect the positivity of $\ct$
with some algebraic structure on the set of partitions itself.
\par
In \l BSp1\r\ we examined one special choice for
$\ct=\hat\ct_\mu$, namely there $\hat\ct_\mu
(\cV)$ was defined with the help of the number of crossings of the
partition $\cV$ and we could reduce the positivity of $\rt$ to a
question on positive definiteness of some function on the
permutation groups $S_r$. The resulting family of Brownian
motions was given by operators $c(f)$ which fulfill the relations
$c(f)c^*(g)-\mu c^*(g)c(f)=<f,g>1$ and give for $\mu$ varying
between -1 and 1 an interpolation between the fermionic and bosonic
relations (comp. also \l BSp2,Spe3,Gre,Fiv,Zag\r).
\par
Here, we shall consider another choice for $\ct$, namely we define it
with the help of the number of connected components of the partitions.
In this case a reduction to the permutation groups $S_r$ is not possible.
\par
In the next section we define the notion of \lq positive definite
function' for 2-partitions and develop some basic theory for such
functions. In particular, we show that this class of functions is
closed under pointwise multiplication.
In Sect. 3 we examine our special example $\ct=\ct_q$, defined with the
help of the number of connected components of the 2-partition and
depending on some parameter $q$. We show the positive definiteness
of $\ct_q$ for $0\leq q\leq 1$ by giving a \lq Fock space' representation
of the corresponding Brownian motion. We also extend, in Sect. 5,
the definition
of $\ct_q$ to negative $q$ and find some connection with our concept
of $\psi$-independence \l BSp3\r. In Sect. 4 we show
that our Brownian
motion is intimately connected with Voiculescu's concepts of
\lq free product' and \lq free convolution'.\par
Since $\ct_1$ gives rise to the bosonic relations and $\ct_{-1}$
corresponds to the fermionic relations we get in this way again an
interpolation between these two cases (including the free case
\l Spe1,Voi1,Maa,KSp\r, which corresponds to $\ct_0$).
\vskip1cm
\heading
{\bf 2. General theory of positive definite functions on 2-partitions}
\endheading
We want to define Brownian motions given by special states with the
help of pairing prescriptions. Our objects of interest consist thus of
pairs $(\cA,\rho)$, where $\cA$ is a unital
$*$-algebra and $\rho$ some special
state on $\cA$.
By a state $\rho$ on a unital $*$-algebra $\cA$ we will always mean
a positive ($\rho(aa^*)\geq 0$ for all $a\in\cA$), hermitean
($\rho(a^*)=\overline{\rho(a)}$ for all $a\in\cA$), and unital
($\rho(1)=1$) linear functional on $\cA$.
In the following $\HH$ will always be an infinite dimensional
separable complex
Hilbert space with an involution $f\mapsto\bar f$. A canonical choice
would be $\HH=L^2(\RR)$ with $\bar f(x)=\overline{f(x)}$.
Now we choose as $\cA$ the free unital $*$-algebra generated by
generators $\omega(f)$ and $\omega(f)^*$ for all $f\in\HH$, divided
by the relations ($\lambda,\mu\in\CC$, $f,g\in\HH$)
 $$\align \omega(\lambda f+\mu g)&=\lambda \omega(f)+\mu \omega(g)\\
\omega(f)^*&=\omega(\bar f),\endalign$$
i.e. $\cA$ is the tensor algebra over $\HH$ (with the canonical
embedding $\omega:\HH\to\cA,\quad f\mapsto \omega(f)$),
made to a $*$-algebra by
putting $\omega(f)^*=\omega(\bar f)$.
We have a natural topology on $\cA\cong\bigoplus_{n\geq 0}\HH^{\otimes n}$
which is given by the (full Fock space) scalar product
 $$<\omega(f_1)\dots \omega(f_n),\omega(g_1)\dots\omega(g_m)>=
\delta_{nm}<f_1,g_1>\dots <f_n,g_n>.$$
In the following we will also use some special subalgebras $\cAn$ of
$\cA$ of the following form:
Choose an orthonormal basis $\{f_i\}$ of $\HH$ with $\bar f_i=f_i$ and
put $\omega_i:=\omega(f_i)$. Then we denote by $\cAn=$ $\koi$ the
unital $*$-subalgebra of $\cA$ generated by
$\omega_i=\omega_i^*$.  By $\{\oo_i\mid i\in\NN\}$ we shall
always denote a set of such generators $\oo_i=\oo(f_i)$, where
$\{f_i\mid i\in\NN\}$ is some orthonormal basis of $\HH$ as above.
Since we are only interested in states $\rho$ on $\cA$ which are
continuous in the natural topology, $\rho$
is determined by its restriction to such a subalgebra $\cAn=\koi$.
\par
For the definition of our special states $\rho$ on $\cA$ we need some
preliminaries about 2-partitions of sets.
Let $S$ be an ordered set. Then we denote by $\Pz(S)$ the set of all
2-partitions of $S$, i.e. $\cV\in\Pz(S)$ if $\cV=\{V_1,\dots,V_r\}$ where
each $V_i$ is an ordered set containing exactly two elements, i.e.
it has the form $V_i=(k,l)$ with $k,l\in S$ and $k<l$, such that all
$V_i$ are disjoint and their union is $S$. In particular, we must have
$\#S=2r$, hence we shall always assume that $S$ has an even number of
elements in order to have $\Pz(S)\not=\emptyset$.
Since only the order of $S$ is important, it is in this respect equivalent
to $(1,2,\dots,\#S)$. We shall use in the following this identification
freely. In particular, functions on $\Pzn$ extend canonically to functions
on $\Pz(S)$ for all ordered $S$ with $\#S=n$. In the next section we shall
also use the {\it set of inversions} $I(\cV)$ of a 2-partition
 $\cV=\{V_1,\dots,V_r\}\in\Pzr$. If we write $V_i=(k_i,l_i)$, then it
is defined as
 $$I(\cV):=\{(i,j)\mid k_i<k_j<l_i<l_j\}.$$
Furthermore, for a 2-partition $\cV\in\Pzr$ we denote by
$\cV^*\in\Pzr$ its \lq adjoint' which is given by reversing the
order of $(1,\dots,2r)$, i.e. for
 $\cV=\{(k_1,l_1),\dots,(k_r,l_r)\}$ we have
 $$\cV^*:=\{(2r+1-l_1,2r+1-k_1),\dots,(2r+1-l_r,2r+1-k_r)\}.$$
\par
Given a
complex-valued function $\ct$ on $\PzN:=\bigcup_{r=1}^\infty \Pzr$,
we define now a state $\rho=\rho_\ct$ on $\cA$ by linear
extension of $\rho_\ct\lb 1\rb=1$ and (for all $n\in\NN$ and
all $f_1,\dots,f_n\in\HH$)
 $$\rho_\ct\lb \oo(f_1)\dots\oo(f_n)\rb=
\cases 0,&\text{if $n$ odd}\\
\sum\Sb\cV=\{V_1,
\dots,V_r\}\\ \in\Pzr\endSb \rho\lb V_1\rb \dots \rho\lb V_r\rb
\cdot \ct(\cV),&\text{if $n=2r$,}\endcases$$
where, for $V=(k,l)$, we put
 $$\rho_\ct\lb V\rb
=\rho\lb V\rb=\rho\lb \oo(f_k)\oo(f_l)\rb =<\bar f_k,f_l>.$$
Note that such states $\rho_\ct$ are continuous on $\cA$.\par
The motivation for this definition comes essentially from the
following general form of a central limit theorem. Another motivation
can be found in \l Spe3\r, where such states appear as models for the
description of statistics of macroscopic fields.
\proclaim{Theorem 0}
Let $\cB$ be an unital $*$-algebra equipped with a state $\ff$.
Consider selfadjoint elements $b_i=b_i^*\in\cB$ ($i\in\NN$) which
fulfill the following assumptions:\newline
i) We have $\ff(b_{i(1)}\dots b_{i(n)})=0$ for all $n\in\NN$ and all
$i(1),\dots,i(n)\in\NN$ with the property that
one of the $i(k)$ is different from all others, i.e. such that
there exists a $k$ ($1\leq k\leq n$) with $i(l)\not=i(k)$ for all
$l\not=k$.\newline
ii) We have an invariance of moments under permutations, i.e. for each
permutation $\pi$ of the natural numbers we have
 $$\ff(b_{i(1)}\dots b_{i(n)})=\ff(b_{\pi(i(1))}\dots b_{\pi(i(n))})$$
for all $n\in\NN$ and all $i(1),\dots,i(n)\in\NN$.\newline
1) Now consider for each $N\in\NN$
 $$S_N:=\frac {b_1+\dots+b_N}{\sqrt N}.$$
Then, for all $n\in\NN$,
 $$\lim_{N\to\infty}\ff(S_N^n)=\cases
0,&\text{if $n$ odd}\\
\sum_{\cV\in\Pz(1,\dots,2r)} \ct(\cV),&\text{if $n=2r$},\endcases$$
where the function $\ct$ is given by the common value
 $\ct(\cV):=\ff(b_{i(1)}\dots b_{i(2r)})$ of all index-tuples
 $(i(1),\dots,i(2r))$ with the property that (with $\cV=\{V_1,\dots,V_r\}$)
 $$i(k)=i(l)\qquad \Longleftrightarrow \qquad \text{there exists a $j$, such
that $k,l\in V_j$}.$$\noindent
2) More generally, consider for all $s,t\in\RR$ with $s<t$ and all $N\in\NN$
 $$S_N(s,t):=\frac 1{\sqrt N} \sum_{i=\lbrack s\cdot N\rbrack +1}^{\lbrack
t\cdot N\rbrack} b_i.$$
Then, for all $n\in\NN$ and all $s_i,t_i\in\RR$ with $s_i<t_i$, we have
 $$\lim_{N\to\infty} \ff(S_N(s_1,t_1)\dots S_N(s_n,t_n))=0$$ for $n$ odd,
and otherwise ($n=2r$)
 $$\multline \lim_{N\to\infty} \ff(S_N(s_1,t_1)\dots S_N(s_n,t_n))=\\
=
\sum_{\cV\in\Pz(1,\dots,2r)} <\chi_{(s_{k_1},t_{k_1})},\chi_{(s_{l_1},
t_{l_1})}>\dots
<\chi_{(s_{k_r},t_{k_r})},\chi_{(s_{l_r},
t_{l_r})}> \ct(\cV),\endmultline$$
where $\cV=\{V_1,\dots,V_r\}$ with $V_i=(k_i,l_i)$ (where $k_i<l_i$)
and $\chi_{(s,t)}\in
L^2(\RR)$ is the characteristic function of the interval $(s,t)\subset\RR$.
\endproclaim
The proof goes along the same lines as the proofs of the central limit
theorems in \l Spe1,Spe2\r, our assumptions guarantee the applicability of
the arguments used there and we refer to these references for details.
This kind of central limit theorem can essentially be traced back to
\l GvW,vWa\r.
\par
One should think of this theorem in the following way: Let, for each $i\in\NN$,
$\cB_i:=\kbi$ be the $*$-algebra generated by $b_i$.
Then $\cB$ should be considered as a \lq product' of the $\cB_i$, which
are lying as \lq independent' subalgebras in $\cB$. Of course, the notion
of \lq independence' is quite subtle and we have to make some comments on
this in the following.\par
Guided by the notion of \lq independence' of K\"ummerer \l K\"um2\r\ we
want to have some factorization properties for our states $\rho_\ct$,
which are
guaranteed if we request the corresponding properties of the
function $\ct$.
There are different possibilities for the notion of independence. The
weakest form for a state $\rho$ on $\cA$ would be the following
factorization:
 $$\rho\lb \oo(f_1)\dots\oo(f_k)\oo(f_{k+1})\dots\oo(f_n)\rb=
\rho\lb\oo(f_1)\dots\oo(f_k)\rb\cdot\rho\lb\oo(f_{k+1})\dots\oo(f_n)\rb$$
for all $k,n\in\NN$ with $k<n$ and all $f_1,\dots,f_n\in\HH$ such that
 $$
<\bar f_i,f_j>=0\qquad \text{for all $i=1,\dots,k$ and $j=k+1,\dots,n$.}$$
A stronger requirement, namely the factorizing of pyramidal
products, is
 $$\multline
\rho\lb\oo(f_1)\dots\oo(f_k)\oo(f_{k+1})\dots\oo(f_l)\oo(f_{l+1})\dots
\oo(f_n)\rb=\\=
\rho\lb\oo(f_{k+1})\dots\oo(f_l)\rb\cdot
\rho\lb\oo(f_1)\dots\oo(f_k)\oo(f_{l+1})\dots\oo(f_n)\rb\endmultline$$
for all $k,l,n\in\NN$ with $k<l<n$ and all $f_1,\dots,f_n\in\HH$ such that
 $$<\bar f_i,f_j>=0\qquad\text{for all $i=1,\dots,k,l+1,\dots,n$ and
$j=k+1,\dots,l$.}$$
The translation of these requirements for $\rt$ to the function $\ct$
is given in the following definition.
\demo{Definition}
i) A function $\ct$ on $\PzN$ is called {\it weakly multiplicative}, if
we have for all $k,n\in\NN$ with $k<n$ and all $\cV_1\in\Pz(1,\dots,k)$
and $\cV_2\in\Pz(k+1,\dots,n)$
 $$\ct(\cV_1\cup\cV_2)=\ct(\cV_1)\cdot\ct(\cV_2).$$
\noindent
ii) A function $\ct$ on $\PzN$ is called {\it multiplicative}, if we
have for all $k,l,n\in\NN$ with $k<l<n$ and all $\cV_1\in\Pz(1,\dots,k,
l+1,\dots,n)$ and $\cV_2\in\Pz(k+1,\dots,l)$
 $$\ct(\cV_1\cup\cV_2)=\ct(\cV_1)\cdot\ct(\cV_2).$$
\enddemo
It is easy to see that $\rt$ has the pyramidal factorization property,
if $\ct$ is multiplicative (and the same for weak factorization property
and weak multiplicativity).
\demo{Example}
For $S=(1,2,3,4)$ we have three 2-partitions, namely
 $$\cV_1=\{(1,2),(3,4)\},\quad \cV_2=\{(1,3),(2,4)\},\quad
\cV_3=\{(1,4),(2,3)\}.$$
Weak multiplicativity gives
 $$\ct(\cV_1)=\ct(1,2)\cdot\ct(3,4)=\ct(1,2)\cdot\ct(1,2),$$
where $\cV=\{(1,2)\}$ is the unique 2-partition of the set $(1,2)$.
Strong multiplicativity gives the additional requirement
 $$\ct(\cV_3)=\ct(1,4)\cdot\ct(2,3)=\ct(1,2)\cdot\ct(1,2)=\ct(\cV_1).$$
For $\cV_2$ there is no equation for a reduction.
\enddemo
We would like to have for our Brownian motions the pyramidal
factorization property,
thus we shall only consider multiplicative $\ct$.
Of course, $\rt$ is only of any use, if it is a state.
\demo{Definition}
A multiplicative
$\ct:\PzN\to\CC$ is called {\it positive definite}, if
$\rho_\ct$, given by
$\rt\lb 1\rb=1$ and
 $$\rt\lb \omega(f_1)\dots\omega(f_n)\rb=\cases
0,&\text{if $n$ odd}\\
\sum\Sb \cV=\{V_1,\dots,V_r\}\\ \in\Pzr\endSb \rho\lb V_1\rb\dots
\rho\lb V_r\rb\cdot\ct(\cV),&\text{if $n=2r$,}\endcases $$
is positive on $\cA$.
\enddemo
Note that positivity of $\rho_\ct$ implies that $\ct$ is
hermitian, i.e.
$\ct(\cV^*)=\overline{\ct(\cV)}.$ \par
The above definition is quite indirect
and one would like to characterize
positive definiteness in terms of some algebraic structure of
$\PzN$. But till now we have not been able to suceed in doing so.\par
Note that because of the continuity of the $\rho_\ct$ it suffices to
have positivity
for the restriction of $\rho_\ct$ to $\cAn=\koi$ for
an arbitrary set of generators $\{\oo_i\mid i\in\NN\}$.\par
We shall now go on to derive some general properties of positive
definite $\ct$, which will show that the assigning of the same name as in
the group case (see, e.g., \l PaS\r)
is not an arbitrary act, but that there are indeed some
similarities. Furthermore, this general theory will also be of some use
in the next section, where we shall consider our special example of
positive definite $\ct$.\par
First of all, one should notice that all permutation groups
$S_r$ are contained in $\PzN$ via
 $$S_r\ni\pi\mapsto \cV_\pi\in\Pzr,$$
where $\cV_\pi:=\{(i,2r+1-\pi(i))\mid i=1,\dots,r\}$.
Thus a function $\ct$ on $\PzN$ gives via restriction also a function
on $S_r$.
\proclaim{Theorem 1}
Let $\ct$ be a positive definite function on $\PzN$. Then,
for all $r\in\NN$, the restriction
of $\ct$ to $S_r$ is also positive definite (in the usual sense).
\endproclaim
\demo{Proof}
Let $h:S_r\to\CC$ be an arbitrary function. Then we have to show
 $$\sum_{\pi,\sigma\in S_r}\ct(\sigma^{-1}\pi)h(\pi)\overline{h(\sigma)}
\geq0.$$
Now let us take some generators $\oo_i=\oo_i^*$ ($i=1,\dots,r$)
and put
 $$a:=\sum_{\pi\in S_r} h(\pi)\oo_{\pi(1)}\dots\oo_{\pi(r)}.$$
The assumed positivity of $\rt$ gives
 $$\rt\lb aa^*\rb=\sum_{\pi,\sigma\in S_r} h(\pi)\overline{h(\sigma)}
\rt\lb\oo_{\pi(1)}\dots\oo_{\pi(r)}\oo_{\sigma(r)}\dots\oo_{\sigma(1)}\rb
\geq0.$$
Since
 $$
\rt\lb\oo_{\pi(1)}\dots\oo_{\pi(r)}\oo_{\sigma(r)}\dots\oo_{\sigma(1)}\rb
=\rt\lb\oo_{\sigma^{-1}\pi(1)}
\dots\oo_{\sigma^{-1}\pi(r)}\oo_{r}\dots\oo_{1}\rb
=\ct(\sigma^{-1}\pi)$$
we get the assertion.\newline
\line{\hfil $\diamondsuit$}
\enddemo
Note that $S_r$ is only a small part of $\Pzr$,
and that in general it
is not sufficient for some function $\ct$ on $\PzN$ to be positive
definite that its restriction to all $S_r$ is positive definite.
We also do not know whether there is some canonical procedure to extend
a collection of
positive definite functions on all $S_r$ to some positive definite
function on $\PzN$.
\demo{Remark}
One might think of the following extension procedure (inspired by our
work in \l BSp1\r): Introduce in the full Fock space $\cF(\HH)$ of $\HH$ a
new scalar product by
 $$<f_1\otimes\dots\otimes f_n,g_1\otimes\dots\otimes g_m>_\ct:=
\delta_{nm}\sum_{\pi\in S_n}<f_1,g_{\pi(1)}>\dots <f_n,g_{\pi(n)}>
\ct(\pi)$$
and write $\oo(f):=c^*(f)+c(\bar f)$, where $c^*(f)=l^*(f)$ is the usual
left creation operator (cf. \l Eva,Voi1\r)
and $c(f)$ is its adjoint with respect to the new
scalar product. If we take now for $\rho$ the vacuum expectation
state then this gives of course some positive definite function on
$\PzN$, but this extension is usually not the \lq right' one, cf. the
corresponding remark for our example in the next section.
\enddemo
Next, we prove an important technical fact, which tells us that the
positive definiteness of $\ct$ depends merely on the values of $\rt$ on
such moments $\oo(f_1)\dots\oo(f_n)$, where only one pairing does
contribute to the evaluation of $\rt$.
\proclaim{Theorem 2}
1) Let $\ct$ be a function on $\PzN$ and assume that
we have for some distinguished set of
generators $\oo_i=\oo_i^*$ ($i\in\NN$) of $\cAn:=\koi\subset\cA$
a state $\rho$ on $\cAn$ which fulfills the following requirements:
\roster
\item"i)"
$\rho\lb \oo_{i(1)}\dots\oo_{i(n)}\rb=0$ for all $n\in\NN$ and all
$i(1),\dots,i(n)\in\NN$ with the property that at least one of the
indices $i(1),\dots,i(n)$ appears exactly once
\item"ii)"
$\rho\lb \oo_{i(1)}
\dots\oo_{i(2r)}\rb=\rho\lb V_1\rb \dots\rho\lb V_r\rb \cdot\ct(\cV)
=\ct(\cV)$
for all $r\in\NN$ and all $i(1),\dots,i(2r)\in\NN$ with the property that
there is exactly one 2-partition $\cV=\{V_1,\dots,V_r\}$ such that
$\rho\lb V_1\rb \dots\rho\lb V_r\rb $ is different from zero, i.e.
each index in $i(1),\dots,i(2r)$ appears exactly twice
\item"iii)"
for all $n\in\NN$ there exists a constant $C_n\in\RR$ such that
 $\vert\rho\lb\oo_{i(1)}\dots\oo_{i(n)}\rb\vert\leq C_n$ for all
$i(1),\dots,i(n)\in\NN$.
\endroster
Then $\ct$ is positive definite.
\newline
2) If we have a state $\rho$ on $\cA$ such that
the properties i), ii), and iii) of part 1) are fulfilled for
$\oo_i=\omega(f_i)$ for all
orthonormal bases $\{f_i=\bar f_i\}$ of $\HH$, then
$\rho=\rt$.
\endproclaim
\demo{Proof}
1) We shall show that we get the general form of $\rt$ via a central
limit theorem if we have the right expression for the special moments
as in our assumption. Note that $\rt$ may be different from $\rho$.
\newline
Let us denote our distinguished set of generators $\oo_i$ by some
identification $\NN\cong\NN\times\NN$ also as $\oo_i^k$ ($i,k\in\NN$).
Then define a new state $\rN$ on $\cAn=\koi$ by
 $$\rN\lb\oo_{i(1)}\dots\oo_{i(n)}\rb:=\rho\lb\oo_{i(1)}^{(N)}\dots
\oo_{i(n)}^{(N)}\rb$$
where
 $$\oo_i^{(N)}:=\frac 1{\sqrt N}\sum_{k=1}^N\oo_i^k.$$
We claim now that $\rN$ tends pointwise to $\rt$ for $N\to\infty$, i.e.
 $$\lim_{N\to\infty}\rN(a)=\rt(a)\qquad\text{for all $a\in\cAn$.}$$
This implies of course that $\rt$ is positive on $\cAn$, thus on $\cA$,
hence that $\ct$ is
positive definite.\newline
Thus it remains to show the convergence of $\rN(a)$. But this is in the
spirit of central limit theorems and our assumptions
are just sufficient to guarantee the applicability of the arguments
used in \l Spe1,Spe2\r. We refer to these references for details.
\newline
2) Consider
$\rho\lb \oo(f_{i(1)})\dots\oo(f_{i(n)})\rb $,
where $\{f_i=\bar f_i\}$ is an
orthonormal basis of $\HH$. Then, for each $N\in\NN$,
we can find another orthonormal basis of $\HH$, enumerated as
$\{f_i^k\mid i\in\NN, k=1,\dots,N\}$, such that
$f_i=\frac 1{\sqrt N}\sum_{k=1}^N f_i^k$.
Inserting these expressions for $f_i$ and letting $N$ tend to $\infty$
we see by the same arguments as in the first part of our proof that
$\rho\lb \oo(f_{i(1)})\dots\oo(f_{i(n)})\rb$
gives exactly the same value as
$\rt\lb \oo(f_{i(1)})\dots\oo(f_{i(n)})\rb $.
\newline
\line{\hfil $\diamondsuit$}
\enddemo
\proclaim{Corollary 1}
If $\ct_1$ and $\ct_2$ are positive definite functions on $\PzN$, then
their pointwise product $\ct$, given by
 $$\ct(\cV):=\ct_1(\cV)\cdot\ct_2(\cV)\qquad (\cV\in\PzN),$$
is positive definite, too.
\endproclaim
\demo{Proof}
Let $\oo_i=\oo_i^*$ be a distinguished set of generators of $\cAn
=\koi$.
According to our assumption $\rho_{\ct_1}$ and $\rho_{\ct_2}$ are states
on $\cAn$. Consider now on $\cAn\otimes\cAn$ the state $\rte\otimes\rtz$ and
embed $\cAn$ in $\cAn\otimes\cAn$ via $\oo_i\mapsto \oo_i\otimes \oo_i$. The
restriction $\rho$ of $\rte\otimes\rtz$ to $\cAn$ is thus a state on $\cAn$
given by linear extension of
 $$\rho\lb \oo_{i(1)}\dots\oo_{i(n)}\rb=
   \rte\lb \oo_{i(1)}\dots\oo_{i(n)}\rb\cdot
  \rtz\lb \oo_{i(1)}\dots\oo_{i(n)}\rb.$$
It fulfills the assumptions of our Theorem 2, because $\rte$ and $\rtz$
fulfill these assumptions. In particular, in case ii) we have
 $$\rho\lb\oo_{i(1)}\dots\oo_{i(2r)}\rb=\ct_1(\cV)\cdot\ct_2(\cV)=\ct(\cV),$$
hence $\rho$ gives rise to the state $\rt$.\newline
Note that $\rho$ depends on the choice of the set of generators and that
this implies $\rho\not=\rt$.\newline
\line{\hfil$\diamondsuit$}
\enddemo
In all known examples of Brownian motions \l HuP,ApH,Par,KSp,BSp1\r\
one can split $\oo(f)$ into
a sum of \lq creation' and \lq annihilation' operators $\oo(f)=c^*(f)+
c(\bar f)$. We can try to imitate this in our general frame. Instead of
$\cA$ we consider now the
free unital $*$-algebra $\cC$ with generators $c(f)$ and
$c^*(f)$ ($f\in\HH$) divided by the canonical linearity relations in
$f$ and $(c(f))^*=c^*(f)$. Again, we can restrict to unital
$*$-subalgebras $\cCn=
\kci$
with generators $c_i$ and $c_i^*$ ($i\in\NN$), where $c_i=c(f_i)$,
$c^*_i=c^*(f_i)$ for some orthonormal basis $\{f_i\}$ of $\HH$ with
$f_i=\bar f_i$ for all $i\in\NN$. Our state $\rt$ should then be replaced
by some state $\rtQh$ on $\cC$ (again determined via
continuity by its restriction to some subalgebra $\cCn$) given by
 $$
\rtQh(c^\sharp(f_1)\dots c^\sharp(f_n))=
\cases 0,&\text{if $n$ odd}\\
\sum\Sb \cV=\{V_1,\dots,V_r\}\\ \in\Pzr\endSb
\rtQh\lb V_1\rb \dots\rtQh\lb V_r\rb
\cdot \ct(\cV),&\text{if $n=2r$,}\endcases
 $$
where now, for $V=(k,l)$, we have
 $$\rtQh\lb V\rb=\rho^Q\lb V\rb=\rho^Q\lb c^\sharp(f_k)c^\sharp
(f_l)\rb=<\bar f_k,f_l> Q(\sharp,\sharp),$$
where $Q$ is the covariance matrix (independent of $i$)
 $$Q=\pmatrix \rho^Q(c_ic_i)&\rho^Q(c_ic_i^*)\\ \rho^Q(c_i^*c_i)&
\rho^Q(c_i^*c_i^*)\endpmatrix.$$
The symbol $\sharp$ denotes the possibility of appearing or not
appearing of a $*$ and an equation with some of these symbols in it
has to be read as a collection of all possible equations where each
of the appearing $\sharp$ is replaced by either $*$ or no $*$, of course
in a consistent way on both sides of the equation.\par
The case $Q=\pmatrix 1&1\\1&1\endpmatrix$ corresponds to our previous
state $\rt$ on $\cA$. We can realize without problem the symmetric
case $Q=\pmatrix 0&1\\1&0\endpmatrix$ on $\cC$ by embedding $\cCn=\kci$
into $\cAn=\koi$
via
 $$c_k=\frac 1{\sqrt 2} (\oo_{2k}+i\oo_{2k+1}),\quad
   c^*_k=\frac 1{\sqrt 2} (\oo_{2k}-i\oo_{2k+1})$$
and take for $\rtQh$ the retract of $\rt$. This may be considered as
a quasi-free state of infinite temperature $T=\infty$.
In analogy with the bosonic
and fermionic cases,
one would like to define also other quasi-free states
as states $\rtQh$ on $\cC$ with more general $Q$.
In particular,
zero temperature $T=0$ would correspond to a Fock or vacuum state
characterized by $Q=\pmatrix 0&1\\0&0\endpmatrix$.
Whether the positive definiteness of $\ct$ implies the positivity of
$\rtQh$ for such states apart from $T=\infty$ is at the
moment an open question.\par
However, in any case we can copy the proofs of Theorem 2 and its
Corollary 1 also for $\rtQh$, i.e. the positivity of $\rtQh$ depends
again merely on the values of $\rtQh$ on special moments as in Theorem 2
and the positivity of $\rho_{\ct_1}^{Q_1}$
and $\rho_{\ct_2}^{Q_2}$
implies the positivity of $\rho_{\ct_1\cdot\ct_2}^{Q_1\circ Q_2}$,
where $Q_1\circ Q_2$ denotes the entrywise (Schur) product of the
covariance matrices $Q_1$ and $Q_2$.
\vskip1cm
\heading
{\bf 3. A special example of a Brownian motion}
\endheading
We shall now make a special choice for the function $\ct$ to get
some new example of a Brownian motion.\par
In \l BSp1\r\ we considered Brownian motions which were given by
special
$\ct=\hat \ct_\mu$ of the form $\hat \ct_\mu(\cV)=\mu^{\#I(\cV)}$, where
$I(\cV)$ is the set of inversions of the 2-partition $\cV$ and $\mu$ is
a parameter. We showed in \l BSp1\r\ (see also \l BSp2,Gre,Fiv,Spe3,Zag\r)
that for $-1\leq\mu\leq 1$ this $\hat \ct_\mu$ is a positive definite
function
on $\PzN$. In this case, the whole problem could be reduced to the question
whether the restriction of $\hat\ct_\mu$ to $S_r$ is positive definite for
all $r\in\NN$.
\par
Here, we shall treat another example. Again, we consider a whole family of
functions $\ct_q$, where $q$ varies now between 0 and 1. In the end we shall
also make some extension of this to negative $q$ and find some connection
to the concept of $\psi$-independence, which was introduced in \l BSp3\r.
Instead of the number of inversions (as for $\hat\ct_\mu$)
we choose now the number of connected
components of a partition $\cV$ for the definition of $\ct_q(\cV)$.
This is in some sense the most canonical form for a multiplicative $\ct$.
If we have $\cV\in\Pzr$, then we say that it decomposes
into $\cV=\cV_1\cup \dots\cup
\cV_k$ (where $\cV_i\in\Pz(S_i)$, with disjoint $S_i$ whose union is
$(1,\dots,2r)$), if we have for all multiplicative $\ct$
the factorization $\ct(\cV)=\ct(\cV_1)\dots\ct(\cV_k)$.
If a partition $\cV_0$ cannot be decomposed into subpartitions, then
we call such a $\cV_0$ {\it connected} or a {\it block}.
If $\cV=\cV_1\cup\dots\cup\cV_k$
is a decomposition of $\cV$ into blocks, then we write
$B(\cV):=\{\cV_1,\dots,\cV_k\}$.
A block is a partition which is connected in its canonical graphic
representation (see \l Spe1\r).
\par
A multiplicative $\ct$ is determined by giving its values on all such
connected blocks. A natural choice for such a function is
 $$\cV_0\mapsto \ct_q(\cV_0)=q^{\#\cV_0-1}\qquad\text{if $\cV_0$ is
a block,}$$
or more general
 $$\ct_q(\cV)=q^{\#\cV-\#B(\cV)}\qquad\text{for $\cV\in\PzN$.}$$
We shall prove that this $\ct_q$ is positive definite for all $q$ with
$0\leq q\leq 1$.\par
\demo{Remark}
For $q=0$ and $q=1$ we get the free and bosonic (classical) Brownian
motions, respectively. Hence, as in \l BSp1\r, we have an interpolation
between the free and bosonic case, but in contrast to \l BSp1\r\ we
cannot include the fermionic case directly in this interpolation. Later,
we shall discuss what can be done for negative $q$.
\enddemo
Let us first describe the restriction of $\ct_q$ to $S_r$. One sees
easily that we can write it in the form $\ct_q(\pi)=q^{d(\pi)}$,
where $d(\pi)$ is the following length function on $S_r$.
Let $e_1,\dots,e_{r-1}$ be the transpositions of neighbouring
elements, i.e. the generators of the Coxeter group $S_r$. Then each
$\pi\in S_r$ can be written (in many ways) as a product of these
generators. In \l BSp1\r\ we used the minimal length $i(\pi)$ of such
representations for $\pi$ (which is equal to the number of inversions
of $\pi$, i.e. to $\#I(\pi):=\#I(\cV_\pi)$)
for the definition of our function $\hat\ct_\mu(\pi)=\mu^{i(\pi)}$.
Here, our $d(\pi)$
is the number of different generators in such a minimal representation.
Although a minimal representation is not unique in general, $d(\pi)$
is well-defined \l Bou\r. For example, in $S_3$ we have $e_1e_2e_1=
e_2e_1e_2$, thus $d(e_1e_2e_1)=2$; whereas $i(e_1e_2e_1)=3$.
Of course, $d(1)=i(1)=0$. From this we see that our functions
$\ct_q$ on $\PzN$ are again extensions of quite
natural functions on $S_r$.
In the case of permutations, $\pi\in S_r$ is
connected or a block if and only if there
is no $k\in\{1,\dots,r-1\}$ such that $\pi(1,\dots,r-k)=(1,\dots,r-k)$.
Thus, with the definition
$B_k:=\{1,2,\dots,k\}$ for $k\in\NN$,
the block structure of a permutation
$\pi\in S_r$ can be described by the set
 $$b(\pi)=\{r-k\mid 1\leq k\leq r-1, \pi(B_{k})=B_{k}\}
\subset\{1,\dots,r-1\}.$$
\demo{Remark}
Let us verify our remark from the last section, that we cannot get
$\ct_q$ from its restrictions to $S_r$ by a mere change of the
scalar product in the full Fock space $$\cF(\HH)=
\bigoplus_{n=0}^\infty \HH^{\otimes n}.$$
One can see directly
(of course, this will also follow from our considerations in Theorem 3),
that the restriction of $\ct_q$ to $S_r$, $\pi\mapsto q^{d(\pi)}$,
is positive definite for all $r\in\NN$, thus
the following definition gives a scalar product on the full Fock space of
$\HH$ ($f_1,\dots,f_n,g_1,\dots,g_m\in\HH$)
 $$<f_1\otimes\dots\otimes f_n,g_1\otimes\dots\otimes g_m>_q:=
\delta_{nm}\sum_{\pi\in S_n}<f_1,g_{\pi(1)}>\dots <f_n,g_{\pi(n)}>
q^{d(\pi)}.$$
Now define, for each $f\in\HH$, $c^*(f)=l^*(f)$ as the usual left
creation operator (cf. \l Eva,Voi1\r), i.e. for $f_1,\dots,f_n\in\HH$
 $$c^*(f)f_1\otimes\dots\otimes f_n=f\otimes f_1\otimes \dots\otimes f_n$$
and take $c(f)$ as the adjoint of $c^*(f)$ with respect to $<\thinspace,
\thinspace>_q$. Then put $\oo(f):=c^*(f)+c(\bar f)$ and define the
state $\rho$ as the vacuum expectation state. This will of course
yield a positive definite function on $\PzN$, let's call it $\tilde \ct_q$,
but it is different from our wanted $\ct_q$, as we can see by
determing their values on
$\cV:=\{(1,4),(2,7),(3,6),(5,8)\}$.
For $\tilde\ct_q$ we have (with $\{f_i\}$ an orthonormal basis of $\HH$
and $c_i=c(f_i)$)
 $$\align
\tilde\ct_q(\cV)&=<\Omega,
\oo_4\oo_2\oo_3\oo_4\oo_1\oo_3\oo_2\oo_1\Omega>_q\\
&=<\Omega,c_4c_2c_3c_4^*c_1c_3^*c_2^*c_1^*\Omega>_q.\endalign$$
Some small calculations give
 $$c(f_1)f_3\otimes f_2\otimes f_1=\frac{q^2}{1+q}(f_3\otimes f_2 +
f_2\otimes f_3),$$
which yields
 $$\align
\tilde\ct_q(\cV)&=<c_4c_3^*c_2^*c_4^*\Omega,c_1c_3^*c_2^*c_1^*\Omega>_q\\
&=\bigl(\frac{q^2}{1+q}\bigr)^2<(f_3\otimes f_2+f_2\otimes f_3),
(f_3\otimes f_2+f_2\otimes f_3)>_q\\
&=\bigl(\frac{q^2}{1+q}\bigr)^2 (2+2q)\\
&=\frac{2q^4}{1+q},\endalign$$
which is different from $\ct_q(\cV)=q^4$.
\enddemo
We shall now prove the positive definiteness of $\ct_q$ -
where $\ct_q(\cV)=q^{\#\cV-\#B(\cV)}$ -
by giving an explicit
construction of operators $\omega(f)$ on some Hilbert space. This
Hilbert space will resemble a Fock space as in our foregoing remark,
but its vectors
will also have to carry some information on connectedness. As usual
we shall split $\omega(f)=c^*(f)+c(\bar f)$ into a sum of creation
and annihilation operators, which must be made to adjoints of each
other by an appropriate choice of the scalar product.\par
As our \lq Fock space' $\cF$ we take the linear combinations of some
distinguished unit vector $\Omega$ (vacuum) and vectors of the form
$(f_1\otimes\dots\otimes f_n,A)$, where $n\in\NN$, $f_i\in\HH$, and
 $A\subset\{1,\dots,n-1\}$; of course, we make the canonical
linear identifications like
 $$\bigl((f_1^{(1)}+f_1^{(2)})\otimes\dots\otimes(f_n^{(1)}+f_n^{(2)})
,A\bigr)
=\sum_{i_1,\dots,i_n=1}^2 (f_1^{(i_1)}\odo f_n^{(i_n)},A).$$
The pair
$(f_1\otimes\dots\otimes f_n,A)$ should be thought of as an $n$-particle
vector, where the particles are grouped into connected blocks. The set
$A$ gives the separation points between these blocks.
We number the separation points from right to left, thus, e.g.,
$(f_1\otimes f_2\otimes f_3,\{1\})$ has the connected blocks
$f_1\otimes f_2$ and $f_3$, whereas $(f_1\otimes f_2\otimes f_3,
\{2\})$ has blocks $f_1$ and $f_2\otimes f_3$.\par
We would like to have a scalar product given by
bilinear extension of
 $$\align
<\Omega,\Omega>_q&=1\\
<\Omega,(\fen,A)>_q&=0\endalign $$ and
 $$\multline
<(\fen,A),(\gem,B)>_q=\\=\delta_{nm}
\sum_{\pi\in S_n}<f_1,g_{\pi(1)}>\dots
<f_n,g_{\pi(n)}> q^{(n-1)-\#\lb A\cap B\cap b(\pi)\rb},\endmultline$$
where $n,m\in\NN$, $f_i,g_j\in\HH$ for $i=1,\dots,n$, $j=1,\dots,m$
and $A\subset\{1,\dots,n-1\}$, $B\subset
\{1,\dots,m-1\}$.\par
Before going on we should check whether the bilinear form $<\thinspace,
\thinspace>_q$ is indeed a scalar product.
\proclaim{Theorem 3}
The bilinear form $<\thinspace,\thinspace>_q$ on $\cF\times\cF$ is, for
$0\leq q\leq 1$, positive.
\endproclaim
\demo{Proof}
Fix $n\in\NN$ and let us denote an element of the form
$f_1\otimes\dots\otimes f_n\in\HH^{\otimes n}$ by $\hat f$. Then we
have to show that for all possible choices of $M\in\NN$, $\hat f_1,\dots
,\hat f_M$ and $A_1,\dots,A_M\subset\{1,\dots,n-1\}$ we have
 $$L:=<\sum_{i=1}^M (\hat f_i,A_i),\sum_{j=1}^M (\hat f_j,A_j)>_q
\quad\geq 0.$$
But this is equal to
 $$\align
L&=\sum_{i,j=1}^M \sum_{\pi\in S_n} <\hat f_i,\pi(\hat f_j)> q^{(n-1)-
\#\lb A_i\cap A_j\cap b(\pi)\rb}\\
&=q^{n-1}\frac 1{n!}\sum_{i,j=1}^M\sum_{\pi,\sigma\in S_n}
<\sigma(\hat f_i),\pi(\hat f_j)> q^{-\#\lb A_i\cap A_j\cap b(\sigma^{-1}\pi)
\rb},\endalign$$
where $\pi(\hat f_j)$ denotes the unitary action of $\pi\in S_n$ on
$\HH^{\otimes n}$, i.e.
 $$\pi(f_1\otimes\dots\otimes f_n)=f_{\pi(1)}\otimes\dots\otimes f_{\pi(n)}.$$
Thus we have to show that the kernel $F$ on $\{1,\dots,M\}\times S_n$, given
by
 $$F\bigl((i,\sigma),(j,\pi)\bigr):=q^{-\#\lb A_i\cap A_j\cap b(\sigma^{-1}
\pi)\rb}$$
(for our fixed choice of $A_1,\dots,A_M$), is positive definite. This
suffices, since the kernel $H$ on $\{1,\dots,M\}\times S_n$, given by
 $$H\bigl((i,\sigma),(j,\pi)\bigr):=<\sigma(\hat f_i),\pi(\hat f_j)>$$
(for our fixed choice of $\hat f_1,\dots, \hat f_M$), is positive definite,
hence the positive definiteness of $F$ implies the one of the pointwise
product $H\cdot F$, which gives at once
 $$\sum_{i,j=1}^M \sum_{\pi,\sigma\in S_n} (H\cdot F)\bigl((i,\sigma),
(j,\pi)\bigr)\geq 0,$$
which we wanted to prove.\newline
So let us show that $F$ is positive definite. This will follow for all
$q$ with $0<q<1$, if we can show that the kernel $\Delta$ on
$\{1,\dots,M\}\times S_n$, given by
 $$\Delta\bigl((i,\sigma),(j,\pi)\bigr):=\#\lb A_i\cap A_j\cap
b(\sigma^{-1}\pi)\rb,$$
is positive definite (see, e.g., \l PaS\r).
Let $\chi_C$ denote the characteristic function of a set $C\subset\NN$ and
introduce on the algebra generated by such functions the positive
definite kernel $<\thinspace,\thinspace>$, which is given by bilinear
extension of
 $$<\chi_C,\chi_D>=\cases 1,& C=D\\
0,& C\not=D.\endcases$$
Then one has $\chi_{b(\sigma^{-1}\pi)}(k)=
<\chi_{\sigma(B_{r-k})},\chi_{\pi(B_{r-k})}
>$, where $B_k=\{1,\dots,k\}$,
and because of $\#C=\sum_{k\in\NN}\chi_C(k)$ we obtain
 $$\align
\Delta\bigl((i,\sigma),(j,\pi)\bigr)&=\sum_{k\in\NN}\chi_{A_i}(k)
\chi_{A_j}(k)\chi_{b(\sigma^{-1}\pi)}(k)\\
&=\sum_{k\in\NN}<\chi_{A_i}(k)\chi_{\sigma(B_{r-k})},\chi_{A_j}(k)
\chi_{\pi(B_{r-k})}>,\endalign$$
which is the sum of positive definite kernels, and hence also positive
definite.
\newline
\line{\hfil$\diamondsuit$}
\enddemo
Now we can define for each $f\in\HH$ a creation operator $c^*(f)$ and
an annihilation operator $c(f)$ by linear extension of
($f,f_i\in\HH$, $A\subset\{1,\dots,n-1\}$)
 $$\align
c^*(f)\Omega&=(f,\emptyset)\\
c^*(f)(\fen,A)&=(f\otimes\fen,A\cup\{n\})\endalign     $$
and
 $$\align
c(f)\Omega&=0\\
c(f)(f_1,\emptyset)&=<f,f_1>\Omega\\
c(f)(\fen,A)&=
\sum_{i=1}^n<f,f_i>(f_1\otimes\dots\otimes \check f_i\otimes
\dots \otimes f_n,A\vert_i)\cdot q^{z(i,A)},\endalign$$
where
 $$\align
z(i,A):&=\cases 0,&\text{if $i=1$ and $n-1\in A$}\\
1,&\text{otherwise,}\endcases\\
A\vert_i:&=\cases A\backslash\{n-1\},&\text{if $i=1$ and $n-1\in A$}\\
A\cap\{1,\dots,n-i\},&\text{otherwise,}\endcases\endalign
 $$
and the symbol $\check f_i$ means that $f_i$ has to be deleted.\par
Let us see whether everything fits nicely and $c(f)$ and $c^*(f)$ are
adjoints of each other. Some care has to be taken since, for $q\not=0$,
they are unbounded operators. But for our algebraic frame-work the
following statement is sufficient.
\proclaim{Theorem 4}
We have for all $\eta,\xi\in\cF$ and all $f\in\HH$
 $$<c^*(f)\eta,\xi>_q=<\eta,c(f)\xi>_q.$$
\endproclaim
\demo{Proof}
It is sufficient to show for all $n\in\NN$, all $f_1,\dots,f_n,g_1,\dots
,g_n\in\HH$, all $A\subset\{1,\dots,n-2\}$, and all $B\subset\{1,\dots,n-1\}$
 $$\multline
<c^*(f_1)(\fzn,A),(\gen,B)>_q=\\=<(\fzn,A),c(f_1)(\gen,B)>_q.\endmultline$$
Let us calculate both sides. The left hand side gives
 $$\align
\text{LHS}&=<(\fen,A\cup\{n-1\}),(\gen,B)>_q\\
&=\sum_{\pi\in S_n}<f_1,g_{\pi(1)}>\dots <f_n,g_{\pi(n)}>
q^{(n-1)-\#\lb (A\cup\{n-1\})\cap B\cap b(\pi)\rb},\endalign$$
whereas the right hand side is equal to
 $$\align
&\text{RHS}=\\
&=\sum_{i=1}^n<f_1,g_i><(\fzn,A),(g_1\otimes\dots\otimes \check
g_i\otimes\dots\otimes g_n,
B\vert_i)>_q q^{z(i,B)}\\
&=\sum_{i=1}^n\sum_{\sigma\in S_{n-1}^{(i)}} <f_1,g_i> <f_2,g_{\sigma(2)}>
\dots <f_n,g_{\sigma(n)}>
q^{(n-2)-\#\lb A\cap B\vert_i\cap b(\sigma)\rb+z(i,B)},\endalign
 $$
where $S_{n-1}^{(i)}$ is the set of all bijections from $\{2,\dots,n\}$ to
$\{1,\dots,\check i,\dots,n\}$ and $b(\sigma)$ is defined by considering
$\sigma$ in a canonical way as an element of $S_{n-1}$. Hence our
assertion follows if we have for all $\pi\in S_n$
 $$\multline
(n-1)-\#\lb(A\cup\{n-1\})\cap B\cap b(\pi)\rb=\\
=(n-2)-\#\lb A\cap B\vert_i\cap b(\sigma)\rb +z(i,B),
\endmultline$$
where $i=\pi(1)$ and $\sigma(j)=\pi(j)$ for $j=2,\dots,n$. But this
follows from the definition of $z(i,B)$ and $B\vert_i$ and the fact that
 $$b(\pi)=\lb b(\sigma)\cup\{n-1\}\rb\cap\{1,\dots,n-i\}.$$
\line{\hfil $\diamondsuit$}
\enddemo
Now we can examine whether these objects give us the right state.
Let $\cC=\cC_q$ be the unital $*$-algebra
generated by all $c(f)$ for $f\in\HH$ and
define on $\cC_q$ the state $\rho_q$ as vacuum expectation
 $$\rho_q(a)=<\Omega,a\Omega>_q=<\Omega,a\Omega>_0\qquad(a\in\cC_q).$$
\proclaim{Theorem 5}
We have for all $n\in\NN$ and all $f_1,\dots,f_n\in\HH$
 $$
\rho_q\lb c^\sharp(f_1)\dots c^\sharp(f_n)\rb =
\cases 0,&\text{if $n$ odd}\\
\sum\Sb \cV=\{V_1,\dots,V_r\}\\ \in\Pzr\endSb
\rho^Q\lb V_1\rb \dots\rho^Q\lb V_r\rb
\cdot \ct_q(\cV),&\text{if $n=2r$,}\endcases
 $$
i.e. $$\rho_q=\rho^Q_{\ct_q}\qquad\text{with covariance matrix}\qquad
Q=\pmatrix0&1\\ 0&0 \endpmatrix.$$
\endproclaim
\demo{Proof}
The vanishing of odd moments follows immediately from the observation
that there must be the same number of creation and annihilation
operators for giving a non-zero vacuum expectation value.\newline
In the other case we use the analogue of the
second part of our Theorem 2 and can thus
restrict to the case where all our $f_i$ are elements of an orthonormal
basis of $\HH$ and
where each $f_i$ appears exactly twice in $\{f_1,\dots,f_{2r}\}$,
hence at most one partition, let's say
$\cV_0$, survives in the sum.
Then $\rho_q\lb c^\sharp(f_1)\dots c^\sharp(f_{2r})\rb $
is only different from
zero if all our pairings connect a $c(f)$ with a $c^*(f)$, i.e. we
first have to create a $f$ before we can annihilate it.
Since
$$c(g)c^*(g)(g_1\odo g_n,A)=<g,g>(g_1\odo g_n,A)$$
for all $A\subset\{1,\dots,n-1\}$ and $g,g_1,\dots,g_n\in\HH$ with
$<g,g_i>=0$ ($i=1,\dots,n$), we have the pyramidal
factorization property for $\rho_q$ and
it suffices to consider the case
where
$\cV_0$ is a block, i.e. $\ct_q(\cV_0)=q^{r-1}$.
In this case $\rho_q\lb c^\sharp(f_1)\dots c^\sharp(f_{2r})\rb $
gives exactly $q^{r-1}$, because,
by the assumption that $\cV_0$ is connected,
each annihilation operator apart from
$c(f_1)$ gives a factor $q$.
Thus the formula is valid in this special case. The
general case follows then as in Theorem 2. (Of course, one can also
check this general case directly, but the writing up is a little bit
cumbersome.)
\newline
\line{\hfil $\diamondsuit$}
\enddemo
\proclaim{Corollary 2}
The function $\ct_q$ on the set $\PzN$ of all 2-partitions, given by
 $$\ct_q(\cV)=q^{\#\cV-\#B(\cV)}\qquad(\cV\in\PzN),$$
is positive definite for $0\leq q\leq1$.
\endproclaim
\vskip1cm
\heading
{\bf 4. Connection with the free product}
\endheading
Next, we shall show that our Brownian motions $(\cC_q,\rho_q)$ are
intimately connected with the reduced free product and free convolution
in the sense of Voiculescu \l Voi1,Voi2\r. Let us recall the relevant
definitions. Assume that $\ff_1$ and $\ff_2$ are states on
unital $*$-algebras $\cB_1$ and $\cB_2$, respectively.
Then consider the algebraic free product $\cB_1\star\cB_2$ (with
identifications of the units) and denote by $j_1$ and $j_2$ the
canonical embeddings of $\cB_1$ and $\cB_2$ into $\cB_1\star\cB_2$,
respectively. One can characterize the (reduced) free product
state $\ff_1\star\ff_2$ on $\cB_1\star\cB_2$ by the following condition:
 $$\ff_1\star\ff_2(j_{l(1)}(a_1)\dots j_{l(n)}(a_n))=0,$$
if $n\in\NN$, $l(1),\dots,l(n)\in\{1,2\}$, $l(1)\not=l(2)\not=\dots\not=
l(n)$, and $a_k\in\cB_{l(k)}$ with
$\ff_{l(k)}(a_k)=0$ for all $k=1,\dots,n$.
If we consider two self-adjoint elements $a_i=a_i^*\in\cB_i$ with
distribution $\nu_i$ with respect to $\ff_i$ ($i=1,2$), i.e.
 $$\ff_i(a_i^n)=\int x^nd\nu_i(x)\qquad (i=1,2),$$
then the distribution $\nu$ of $j_1(a_1)+j_2(a_2)$ with respect
to $\ff_1\star\ff_2$ depends only on $\nu_1$ and $\nu_2$ and
is called the free convolution of $\nu_1$ and $\nu_2$, denoted by
$\nu=\nu_1\sqp\nu_2$ \l Voi1,Voi2,Maa,BV\r.
\proclaim{Theorem 6}
Let $\{c_i,c_i^*\mid i\in\NN\}$ denote a distinguished set of generators
of the unital $*$-algebra $\cCn=\kci$
and let $q_1,q_2$ be real numbers with $0\leq q_1,q_2\leq 1$.
Embed $\cCn$ in $\cCn\star\cCn$ via ($i\in\NN$)
 $$c_i\mapsto \sqrt{\frac q{q_1}} j_1(c_i)+\sqrt{\frac q{q_2}}j_2(c_i),\quad
c_i^*\mapsto \sqrt{\frac q{q_1}} j_1(c_i^*)+\sqrt{\frac q{q_2}}j_2(c_i^*)$$
and let $\rho$ be the restriction of $\rqe\star\rqz$ to $\cCn$. Then
 $$\rho=\rho_q,\qquad\text{where}\quad \frac 1q=\frac 1{q_1}+\frac 1{q_2}.
 $$
\endproclaim
\demo{Proof}
Define $q$ by the equation
$\frac 1q=\frac 1{q_1}+\frac 1{q_2}$. Then we have to show for all $n\in\NN$
and all $i(1),\dots,i(n)\in\NN$
 $$\multline
\rqe\star\rqz\bigl\lb\bigl(\sqrt{\frac q{q_1}}j_1(c_{i(1)}^\sharp)+
                 \sqrt{\frac q{q_2}}j_2(c_{i(1)}^\sharp)\bigr)\dots
\bigl(\sqrt{\frac q{q_1}}j_1(c_{i(n)}^\sharp)+
                 \sqrt{\frac q{q_2}}j_2(c_{i(n)}^\sharp)\bigr)\bigr\rb=\\
\quad\\=\rho_q\lb c_{i(1)}^\sharp\dots c_{i(n)}^\sharp\rb.\endmultline$$
For odd $n$ both sides are zero, so we may restrict to $n=2r$.
For the proof we have to use the machinery of \lq non-crossing
cumulants', which was introduced in \l Spe4\r, see also \l NSp\r.
We give here only a short sketch of the main ideas, for the special
definitions and more details we refer to \l Spe4,NSp\r. For a state
$\ff$ on a unital $*$-algebra $\cB$ we consider quantities
$k(\ff)\lb a_1,\dots,a_n\rb$ (for all $a_1,\dots,a_n\in\cB$),
called non-crossing cumulants, which are determined by the
moments $\ff\lb a_1\dots a_n\rb$ ($a_i\in\cB$) of $\ff$ via the
relation
 $$\ff\lb a_1\dots a_n\rb=\sum\Sb \cV=\{V_1,\dots,V_p\}\\ \in
NC(1,\dots,n)\endSb k(\ff)\lb a_{V_1}\rb\dots k(\ff)\lb a_{V_p}\rb.$$
In this formula
the sum runs over all non-crossing partitions $\cV$ of $(1,\dots,n)$ and
for $V=(v_1,\dots,v_s)$ ($v_1<\dots <v_s$) we denote
$k(\ff)\lb a_V\rb:=k(\ff)\lb a_{v_1},\dots,a_{v_s}\rb$.
These
non-crossing cumulants have the crucial property that they linearize
free convolution, i.e. the cumulant of $\ff_1\star\ff_2$
is given by the \lq direct sum' of the
cumulants of $\ff_1$ and of $\ff_2$,
which means
 $$\multline
k(\ff_1\star\ff_2)\lb j_{l(1)}(a_1),\dots, j_{l(n)}(a_n
)\rb=\\
=\cases
k(\ff_1)\lb a_1,\dots,a_n\rb,&\text{if
$l(1)=\dots=l(n)=1$}\\
k(\ff_2)\lb a_1,\dots,a_n\rb,&\text{if
$l(1)=\dots=l(n)=2$}\\
0,&\text{otherwise}\endcases\endmultline$$
for all $l(1),\dots,l(n)\in\{1,2\}$ and all $a_k\in \cB_{l(k)}$.
It is easy to check that in the case where the moments of $\ff$ are
given by a formula involving
summation over all 2-partitions, the non-crossing cumulants
are given by summation over all connected 2-partitions, i.e.
in our case $\cB=\cC_q$ and $\ff=\rho_q$ the non-crossing
cumulant $k_q:=k(\rho_q)$ calculates as
$$\align
k_q\lb c^\sharp_{i(1)},\dots, c^\sharp_{i(2r)}\rb&=\sum\Sb
\cV_0=\{V_1,\dots,V_r\}\\ \in\Pzr\\ \text{$\cV_0$ connected}\endSb
\rho^Q\lb V_1\rb\dots\rho^Q\lb V_r\rb\cdot\ct_q(\cV_0)\\&=
\sum\Sb
\cV_0=\{V_1,\dots,V_r\}\\ \in\Pzr\\ \text{$\cV_0$ connected}\endSb
\rho^Q\lb V_1\rb\dots\rho^Q\lb V_r\rb\cdot q^{r-1}\endalign$$
But this implies that the cumulant
$k(\rqe\star\rqz)$ of the expression
 $$\Bigl\lb\bigl(\sqrt{q/{q_1}}j_1(c_{i(1)}^\sharp)+
                 \sqrt{q/{q_2}}j_2(c_{i(1)}^\sharp)\bigr),\dots,
\bigl(\sqrt{q/{q_1}}j_1(c_{i(2r)}^\sharp)+
                 \sqrt{q/{q_2}}j_2(c_{i(2r)}^\sharp)\bigr)\Bigr\rb$$
is equal to
 $$\align
(q/{q_1})^r& k_{q_1}\lb c^\sharp_{i(1)},\dots, c^\sharp_{i(2r)}\rb+
(q/{q_2})^r k_{q_2}\lb c^\sharp_{i(1)},\dots, c^\sharp_{i(2r)}\rb\\
&=\sum_{\text{connected $\cV_0$}}\rho^Q\lb V_1\rb\dots\rho^Q\lb V_r\rb
\cdot \bigl\{
(q/{q_1})^rq_1^{r-1}+(q/{q_2})^rq_2^{r-1}\bigr\}\\
&=\sum_{\text{connected $\cV_0$}}
\rho^Q\lb V_1\rb\dots\rho^Q\lb V_r\rb\cdot \bigl\{q^r(\frac
1{q_1}+\frac 1{q_2})\bigr\}\\
&=\sum_{\text{connected $\cV_0$}}
\rho^Q\lb V_1\rb\dots\rho^Q\lb V_r\rb\cdot q^{r-1}\\
&=k_q\lb c^\sharp_{i(1)},\dots, c^\sharp_{i(2r)}\rb.
\endalign$$
Thus we have equality of the cumulants, which implies equality of the
states, since moments and cumulants determine each other (see
\l Spe4\r).
\newline
\line{\hfil $\diamondsuit$}
\enddemo
For $q=1/N$ ($N\in\NN$), this gives us the following realization of
our Brownian motion $(\cC,\rho_{1/N})$ with the help of bosonic
operators.
\proclaim{Corollary 3}
Let $a_i^{(k)}$ ($i\in\NN$, $k=1,\dots,N$) for fixed $k$
be generators of the bosonic
relations which are \lq freely independent' for $k\not= l$, i.e.
for all $i,j\in\NN$ and $k,l=1,\dots,N$ with $k\not= l$ we have
 $$\align
a_i^{(k)}a_j^{(k)}-a_j^{(k)}a_i^{(k)}&=0\\
a_i^{(k)}a_j^{(k)*}-a_j^{(k)*}a_i^{(k)}&=\delta_{ij}\endalign$$
and
 $$a_i^{(k)}a_j^{(l)*}=0,$$
and let $\rho$ be the vacuum expectation state on the $*$-algebra
generated by all $a_i^{(k)}$, i.e. $\rho(a)=<\Omega,a\Omega>$, where
$\Omega$ is the vacuum vector characterized by $a_i^{(k)}\Omega =0$
for all $i\in\NN$ and $k=1,\dots,N$. Now embed $\cCn=\kci$
in this algebra as
 $$c_i:=\frac 1{\sqrt N}\sum_{k=1}^N a_i^{(k)},\qquad
   c_i^*=\frac 1{\sqrt N}\sum_{k=1}^N a_i^{(k)*}.$$
Then the restriction of $\rho$ to $\cCn$ gives $\rho_{1/N}$.
\endproclaim
We can also specialize our theorem to an assertion about the \lq Gaussian
measures' connected to our Brownian motions. Let us denote
the spectral measure of $c_i^*+c_i$ with respect to $\rho_q$,
which is of course independent of $i$,
by
$\mu_q$. Then we know that $\mu_0$ is the Wigner semicircle
 $$d\mu_0(x)=\frac 1\pi\sqrt{1-(x/2)^2}dx\qquad x\in \lb -2,2\rb$$
and $\mu_1$ the
usual Gaussian measure
 $$d\mu_1(x)=\frac 1{\sqrt{2\pi}}\exp(-x^2/2)dx\qquad x\in\RR.$$
If we denote by $D_\lambda$ the dilation of probability
measures on $\RR$ by a factor $\lambda$, i.e.
 $$(D_\lambda\mu)(A)=\mu(\lambda^{-1}A) \qquad\text{for $A\subset\RR$
measurable,}$$
and by $\sqp$ the free convolution according to Voiculescu \l Voi2\r,
then we
have the following corollary of our theorem.
\proclaim{Corollary 4}
For all $0\leq q,q_1,q_2\leq 1$ with $\frac 1q=\frac 1{q_1}
+\frac 1{q_2}$ we have
 $$\mu_q=D_{\sqrt{q/q_1}}\thinspace\mu_{q_1}\sqp
         D_{\sqrt{q/q_2}}\thinspace\mu_{q_2}.$$
In particular,
 $$\mu_{1/N}=D_{\sqrt{1/N}}\thinspace\mu_1\sqp\dots\sqp
             D_{\sqrt{1/N}}\thinspace\mu_1\qquad\text{($N$ summands)}$$
for all $N\in\NN$.
\endproclaim
The measures $\mu_{1/N}$ converge for $N\to\infty$ to the Wigner measure
$\mu_0$ and $\mu_{1/N}$ are the (integer) steps in a
free central limit theorem starting with a Gaussian distribution. Thus our
Brownian motion gives a canonical interpolation for the \lq non-integer'
steps of this procedure.
\par
\vskip1cm
\heading
{\bf 5. Extension to negative $q$}
\endheading
Let us finally discuss, whether we can extend our Brownian motion to
negative $q$. One can check that the restriction of our $\ct_q$ to
$S_r$ is not positive definite for negative $q$ as $r\to\infty$.
 Thus a direct
extension is not possible. But we can define $\ct_q$ for negative $q$ by
 $$\ct_q:=\ct_{-q}\cdot\ct_{-1}\qquad (-1\leq q\leq 0),$$
where
 $$\ct_{-1}(\cV):=(-1)^{\#I(\cV)}.$$
This $\ct_{-1}=\hat\ct_{-1}$
gives rise to the fermionic relations and hence is positive
definite \l BSp1\r.
Thus Corollary 1 ensures that the so defined $\ct_q$ is positive
definite, too. \par
One sees easily that Theorem 6 remains true with this definition also
for $-1\leq q,q_1,q_2\leq0$.
The multiplication with $\ct_{-1}$ has the effect that the
bosonic relations ($q=1$) are replaced by the fermionic ones ($q=-1$).
Hence our corollaries are replaced in the following way.
\proclaim{Corollary 5}
Let $b_i^{(k)}$ ($i\in\NN$, $k=1,\dots,N$) for fixed $k$
be generators of the fermionic
relations which are \lq freely independent' for $k\not= l$, i.e.
for all $i,j\in\NN$ and $k,l=1,\dots,N$ with $k\not= l$ we have
 $$\align
b_i^{(k)}b_j^{(k)}+b_j^{(k)}b_i^{(k)}&=0\\
b_i^{(k)}b_j^{(k)*}+b_j^{(k)*}b_i^{(k)}&=\delta_{ij}\endalign$$
and
 $$b_i^{(k)}b_j^{(l)*}=0$$
and let $\rho$ be the vacuum expectation state on the $*$-algebra
generated by all $b_i^{(k)}$, i.e. $\rho(b)=<\Omega,b\Omega>$, where
$\Omega$ is the vacuum vector characterized by $b_i^{(k)}\Omega =0$
for all $i\in\NN$ and $k=1,\dots,N$. Now embed $\cCn=\kci$
in this algebra as
 $$c_i:=\frac 1{\sqrt N}\sum_{k=1}^N b_i^{(k)},\qquad
   c_i^*=\frac 1{\sqrt N}\sum_{k=1}^N b_i^{(k)*}.$$
Then the restriction of $\rho$ to $\cCn$ gives $\rho_{-1/N}$.
\endproclaim
\proclaim{Corollary 6}
We have
 $$\mu_{-1/N}=D_{\sqrt{1/N}}\thinspace\mu_{-1}\sqp\dots\sqp
             D_{\sqrt{1/N}}\thinspace\mu_{-1}\qquad\text{($N$ summands)}$$
for all $N\in\NN$.
\endproclaim
One should note, that $\mu_{-1}$ is nothing else than $\frac 12
(\delta_{-1}+\delta_{+1})$. Thus our $\mu_{-1/N}$ are \lq free
binomial distributions' and $\mu_q$ for negative $q$ is again a
canonical interpolation between the integer steps of the corresponding
free Moivre-Laplace central limit theorem.
\par
Interestingly, in the case of negative $q$
there is also some connection with the notion of
$\psi$-independence.
This is a generalization of the free product \l Boz,BSp3\r: If we
have pairs of states $(\ff_1,\psi_1)$ and $(\ff_2,\psi_2)$ on
unital $*$-algebras $\cB_1$ and $\cB_2$, respectively,
then the state $\ff:=(\ff_1,\psi_1)\star
(\ff_2,\psi_2)$ on $\cB_1\star\cB_2$ is characterized by the condition:
 $$\ff(j_{l(1)}(a_1)\dots j_{l(n)}(a_n))=\ff_{l(1)}(a_1)\dots
\ff_{l(n)}(a_n)$$
if $n\in\NN$, $l(1),\dots,l(n)\in\{1,2\}$, $l(1)\not= l(2)\not=\dots
\not=l(n)$ and $a_k\in\cB_{l(k)}$ with $\psi_{l(k)}(a_k)=0$ for all
$k=1,\dots,n$. Again, we have denoted by $j_1$ and $j_2$ the canonical
embeddings of $\cB_1$ and $\cB_2$ into $\cB_1\star\cB_2$.\par
In \l BSp3,BLS\r\ we calculated the distributions
$\mu_{\alpha,\beta}$ appearing in the central limit
theorem for $\psi$-independence. The moments of $\mu_{\alpha,\beta}$
are given by
 $$\int_\RR x^n d\mu_{\alpha,\beta}(x)=\cases
0,&\text{if $n$ odd}\\
\sum\Sb \cV=\{V_1,\dots,V_r\}\\ \in \Pza\endSb \alpha^{
2\cdot\#B^{(out)}(\cV)}
\beta^{2\cdot\#B^{(in)}(\cV)},&\text{if $n=2r$,}\endcases $$
where $\Pza$ denotes special 2-partitions, namely
\lq non-crossing' partitions which were introduced by Kreweras
\l Kre\r.
In \l Spe1\r\  we called them \lq admissible' partitions.
They may be defined by the
requirement
 $$\cV\in\Pza\Longleftrightarrow \#B(\cV)=r
\Longleftrightarrow B(\cV)=\cV.$$
The blocks $B(\cV)=\cV$ of such partitions are divided into outer and
inner blocks, $B(\cV)=B^{(out)}(\cV)\cup B^{(in)}(\cV)$. A block
$V_j=(k,l)\in\cV$ is called {\it inner}, if there exists another block
$V_r=(\hat k,\hat l)\in\cV$ such that $\hat k<k<l<\hat l$. Otherwise the
block $V_j$ is called {\it outer}.
\par
Now we have the following coincidence.
\proclaim{Theorem 7}
We have
 $$\mu_q=\mu_{1,\sqrt{1+q}}\qquad\text{for $-1\leq q\leq0$.}$$
In particular,
 $$\mu_q=\alpha_q\thinspace(\delta_{\sqrt{1/{\vert q\vert}}}+
                  \delta_{-\sqrt{1/{\vert q\vert}}})+\tilde \mu_q$$
with
 $$d\tilde\mu_q(t)=\chi_{\lb-2\sqrt{1+q}
,2\sqrt{1+q}\rb}(t)\cdot\frac 1{2\pi}
\frac{\sqrt{4(1+q)-t^2}}{1-\vert q\vert t^2}dt$$
and
 $$\alpha_q=\frac 1{4\vert q\vert} \max(1-2(1+q),0).$$
\endproclaim
\demo{Proof}
Let $-1\leq q\leq 0$. Denote again by $\rho_q$ the state
corresponding to $\ct_q$. We have a similar representation on a Fock space
as for positive $q$, one only has to take into account the minus sign in
the definition of the annihilation operator:
 $$\multline
c(f)(\fen,A)=\\
=\sum_{i=1}^n<f,f_i>(f_1\otimes\dots\otimes \check f_i\otimes
\dots f_n,A\vert_i)\cdot (-q)^{z(i,A)}\cdot(-1)^{i-1},
\endmultline$$
which results of course also in a change of the scalar product:
 $$\multline
<(\fen,A),(\gem,B)>_q=\\=\delta_{nm}
\sum_{\pi\in S_n}<f_1,g_{\pi(1)}>\dots
<f_n,g_{\pi(n)}> (-q)^{(n-1)-\#\lb A
\cap B\cap b(\pi)\rb}\cdot(-1)^{\#I(\pi)}.\endmultline$$
Let $c:=c(f)$ for some fixed $f\in\HH$ with $\Vert f\Vert=1$. Consider
now a product $c^\sharp\dots c^\sharp$ of length $n$. For odd $n$,
$\rqcc$ is clearly zero. Thus it is sufficient to show that for all such
products of length $2r$
 $$\rqcc=\sum\Sb\cV=\{V_1,\dots,V_r\}\\
  \in\Pza\endSb \rho^Q\lb V_1\rb\dots
\rho^Q\lb V_r\rb\cdot (1+q)^{B^{(in)}(\cV)}.$$
If $\rqcc$ is different from zero, then there is exactly one
non-crossing partition,
let's say $\cV_0$, such that $\rho^Q\lb V_1\rb\dots\rho^Q
\lb V_r\rb$ does not vanish, i.e. the sum reduces to one summand, and
we are done if we can show that
 $$\rqcc=<\Omega,c^\sharp\dots c^\sharp\Omega>_q
=(1+q)^{\#B^{(in)}(\cV_0)}.$$
That this is indeed true, follows from the following observation:
Let $A,B\subset\{1,\dots,n-1\}$. Then the scalar product
 $$
<(f\otimes\dots\otimes f,A),(f\otimes\dots\otimes f,B)>_q
=\sum_{\pi\in S_n}(-q)^{(n-1)-\#\lb A\cap B\cap b(\pi)\rb}\cdot
(-1)^{\#I(\pi)}$$
vanishes if at least one of the two $n$-particle vectors contains a
connected block, i.e. if $A\cap B\not=\{1,\dots,n-1\}$, because then we
can collect the permutations in $S_n$ in pairs $(\pi,\sigma)$, such that
the contributions of $\pi$ and $\sigma$ in the above sum cancel each other:
If $k\not\in A\cap B$, we pair $\pi$ with $\sigma$ if
 $$\align
\pi(i)&=\sigma(i)\qquad\qquad\text{for $i\not=n-k,n-k+1$}\\
\pi(n-k)&=\sigma(n-k+1)\\
\pi(n-k+1)&=\sigma(n-k),\endalign$$
which implies
 $$\#\lb A\cap B\cap b(\pi)\rb=\#\lb A\cap B\cap b(\sigma)\rb$$
and
 $$\#I(\pi)=\#I(\sigma)\pm 1.$$
{}From the vanishing of scalar products of connected vectors it follows that
for the calculation of
$<\Omega,c^\sharp\dots c^\sharp\Omega>_q$ we only have to
take care of vectors which are not connected, which means that we only have
to pay attention to
the annihilation of the first and the second factor in a tensor product. Thus,
if we act with $c(f)$ on a
(not connected) vector in the $n$-particle space with $n\geq 2$
we get a factor $(1+q)$, an annihilation of a vector in the $1$-particle
space gives a factor 1. Since we have to annihilate $n$-particle vectors
with $n\geq 2$ exactly as often as we have inner blocks of $\cV_0$ we
get the wanted equality.\newline
The explicit formula for $\mu_q$ follows then from \l BSp3,BLS\r.
\newline
\line{\hfil $\diamondsuit$}
\enddemo
Note that we have proved with the help of our concrete representation
the following combinatorial fact.
\proclaim{Corollary 7}
We have for all $r\in\NN$ and all $q$ with $0\leq q\leq 1$
 $$\sum_{\cV\in\Pzr}q^{\#\cV-\#B(\cV)}\cdot(-1)^{\#I(\cV)}=
\sum_{\cV\in\Pza}(1-q)^{\#B^{(in)}(\cV)}.$$
\endproclaim
\demo{Examples}
1) Let us again
consider the three 2-partitions for $2r=4$ as given in Sect. 2.
Thus
 $$\Pz(1,2,3,4)=\{\cV_1,\cV_2,\cV_3\},
\qquad \Pzaa(1,2,3,4)=\{\cV_1,\cV_3\}.$$
The left and right hand side of the corollary are then
 $$\align
LHS&=\ct_{-q}(\cV_1)+\ct_{-q}(\cV_2)+\ct_{-q}(\cV_3)=1-q+1=2-q\\
RHS&=(1-q)^{\#B^{(in)}(\cV_1)}+(1-q)^{\#B^{(in)}(\cV_3)}=(1-q)^0+(1-q)^1
=2-q.\endalign$$\noindent
2) The foregoing example is somehow misleading since usually there are
cancellations on the left hand side of our equation, e.g., for $2r=6$ we
get
 $$RHS=1+2(1-q)+2(1-q)^2=5-6q+2q^2,$$
i.e. on the left hand side two of the 15 summands must cancel each other.
\enddemo
\vskip1cm
\heading
{\bf Acknowledgements}
\endheading
This work has been supported by
Polish National Grant, KBN 4019 (M.B.) and by
the Deutsche Forschungsgemeinschaft,
SFB 123 (R.S.). R.S. thanks Ken Dykema for interesting discussions, by which
Theorems 6 and 7 of this work were inspired. We would also like to thank
an
anonymous referee for pointing out some ambiguities and incorectness
in an earlier version of the paper.
\vskip1cm
\heading
{\bf References}
\endheading
\roster
\item"\l AFQ\r"
L. Accardi, F. Fagnola, and J. Quaegebeur, A representation free quantum
stochastic calculus, {\it J. Funct. Anal.} {\bf 104} (1992), 149-197.
\item"\l ApH\r"
D.B. Applebaum and R.L. Hudson, Fermion Ito's formula and stochastic
evolutions, {\it Comm. Math. Phys} {\bf 96} (1984), 473-496.
\item"\l BV\r"
H. Bervovici and D. Voiculescu, Free Convolution of Measures with
Un\-bounded Support, preprint PAM 572, Berkeley 1992.
\item"\l Bou\r"
N. Bourbaki, {\it Groupes et algebres de lie, Chap. 4,5,6}, Hermann 1968.
\item"\l BLS\r"
M. Bo\D zejko, M. Leinert and R. Speicher, Convolution and limit
theorems for conditionally free random variables,
preprint, Heidelberg 1994.
\item"\l BSp1\r"
M. Bo\D zejko and R. Speicher, An Example of a Generalized Brownian
Motion, {\it Comm. Math. Phys.} {\bf 137} (1991), 519-531.
\item"\l BSp2\r"
M. Bo\D zejko and R. Speicher, An example of a generalized Brownian
motion II, in {\it Quantum Probability and Related Topics VII},
World Scientific, Singapore 1992, 67-77.
\item"\l BSp3\r"
M. Bo\D zejko and R. Speicher, $\psi$-Independent and Symmetrized White
\newline Noises, in {\it Quantum Probability and Related Topics VI},
World Scientific, Singapore 1991, 219-236.
\item"\l Boz\r"
M. Bo\D zejko, Positive definite functions on the free group and the
non commutative Riesz product, {\it Boll. Un. Math. Ital.} {\bf 5A}
(1986),13-22.
\item"\l Eva\r"
D.E. Evans, On $O_n$, {\it Publ. RIMS} {\bf 16} (1980), 915-927.
\item"\l Fiv\r"
D. Fivel, Interpolation between Fermi and Bose statistics using
generalized commutators, {\it Phys. Rev. Lett.} {\bf 65} (1990), 3361-3364.
\item"\l GvW\r"
N. Giri and W. von Waldenfels, An algebraic version of the
central limit theorem, {\it Z. Wahrscheinlichkeitstheorie verw. Gebiete}
{\bf 42} (1978),
129-134.
\item"\l Gre\r"
O.W. Greenberg, Particles with small violations of Fermi or Bose statistics,
{\it Phys. Rev. D} {\bf 43} (1991), 4111-4120.
\item"\l HuP\r"
R.L. Hudson and K.R. Parthasarathy, Quantum Ito's formula and stochastic
evolution, {\it Comm. Math. Phys} {\bf 93} (1984), 301-323.
\item"\l Kre\r"
G. Kreweras, Sur les partitions non croisees d'un cycle, {\it Discr. Math.}
{\bf 1} (1972),
333-350.
\item"\l K\"um1\r"
B. K\"ummerer, Markov dilations on $W^*$-algebras, {\it J. Funct. Anal.}
{\bf 63} (1985),
139-177.
\item"\l K\"um2\r"
B. K\"ummerer,
Markov dilations and non-commutative Poisson processes,
preprint.
\item"\l KPr\r"
B. K\"ummerer and J. Prin, Generalized white noise and non-commutative
stochastic integration, preprint.
\item"\l KSp\r"
B. K\" ummerer and R. Speicher, Stochastic Integration on the Cuntz Algebra
$O_\infty$, {\it J. Funct. Anal.} {\bf 103} (1992), 372-408.
\item"\l LPo\r"
R. Lenczewski and K. Podgorski, A $q$-Analog of the Quantum Central Limit
Theorem for $SU_q(2)$, {\it J. Math. Phys.} {\bf 33} (1992), 2768-2778.
\item"\l Maa\r"
H. Maassen, Addition of Freely Independent Random Variables,
{\it J. Funct. Anal.}
{\bf 106} (1992), 409-438.
\item"\l NSp\r"
P. Neu and R. Speicher, A self-consistent master equation and a new
kind of cumulants, {\it Z. Phys. B} {\bf 92} (1993), 399-407.
\item"\l Par\r"
K.R. Parthasarathy, {\it An introduction to quantum stochastic calculus},
\newline Birkh\"auser 1992.
\item"\l PaS\r"
K.R. Parthasarathy and K. Schmidt, {\it Positive Definite Kernels, Continuous
Tensor Products, and Central Limit Theorems of Probability Theory},\newline
Springer-Verlag (LNM 272), Heidelberg 1972.
\item"\l Spe1\r"
R. Speicher, A New Example of \lq Independence' and \lq White Noise',
{\it Probab. Th. Rel. Fields} {\bf 84} (1990), 141-159.
\item"\l Spe2\r"
R. Speicher, A non-commutative central limit theorem, {\it Math. Z.}
{\bf 209} (1992),
55-66.
\item"\l Spe3\r"
R. Speicher, Generalized Statistics of Macroscopic Fields,
{\it Lett. Math. Phys.} {\bf 27} (1993), 97-104.
\item"\l Spe4\r"
R. Speicher, Multiplicative functions on the lattice of
non-crossing partitions and free convolution,
{\it Math. Ann.} {\bf 298} (1994), 611-628
\item"\l Voi1\r"
D. Voiculescu, Symmetries of some reduced  free product
$C^*$-algebras, in
{\it Operator Algebras and their Connection with
Topology and Ergodic Theory},
LNM 1132, Springer, Heidelberg 1985, 556-588.
\item"\l Voi2\r"
D. Voiculescu, Addition of certain non-commuting random
variables, {\it J. Funct. Anal.} {\bf 66} (1986), 323-346.
\item"\l vWa\r"
W. von Waldenfels, An algebraic central limit theorem in the
anti-com\-muting case, {\it Z. Wahrscheinlichkeitstheorie verw. Gebiete}
{\bf 42} (1978),
135-140.
\item"\l Zag\r"
D. Zagier, Realizability of a Model in Infinite Statistics, {\it Comm. Math.
Phys.} {\bf 147} (1992), 199-210.
\endroster
\enddocument